\documentclass{rspublic}

\usepackage{amsmath}
\usepackage{amssymb}

\newcommand{\Bbf}{\mathbb}

\newcommand{\BB}{{\cal B}}

\newcommand{\FF}{{\cal F}}

\newcommand{\MM}{{\cal M}}

\newcommand{\PP}{{\cal P}}

\newcommand{\concat}{\kern-.25pt\raise4pt\hbox{$\frown$}\kern-.25pt}

\newcommand{\sub}{\subseteq}
\newcommand{\alg}{\mathfrak A}

\newcount\skewfactor
\def\mathunderaccent#1#2 {\let\theaccent#1\skewfactor#2
\mathpalette\putaccentunder}
\def\putaccentunder#1#2{\oalign{$#1#2$\crcr\hidewidth
\vbox to.2ex{\hbox{$#1\skew\skewfactor\theaccent{}$}\vss}\hidewidth}}

\newcommand{\eps}{\varepsilon}

\newcommand{\algB}{\mathfrak B}
\newcommand{\algb}{\mathfrak B}

\newcommand{\intn}{\mbox{\rm int}}

\begin{document}

\title{Measure Recognition Problem}

\author{Mirna D\v{z}amonja}
\affiliation{School of Mathematics, University of East Anglia, Norwich NR47TJ, UK, h020@uea.ac.uk} 


\maketitle

\begin{abstract}{set theory, Boolean algebras, measure, combinatorial characterization}
This is an article in mathematics, specifically in set theory. 
On the example of the Measure Recognition
Problem (MRP) the article
highlights the phenomenon of the utility of a multidisciplinary mathematical approach to
a single mathematical problem, in particular the value of a set-theoretic analysis.
MRP asks if for a given Boolean algebra $\algB$ and a property
$\Phi$ of measures one can recognize by purely combinatorial means
if $\algB$ supports a strictly
positive measure with property $\Phi$. The most famous instance of this problem
is MRP(countable additivity), and in the first part of the article we survey the
known results on this and some other problems. We show how these results naturally
lead to asking about two other specific instances of the problem MRP, namely
MRP(nonatomic) and MRP(separable). Then we show how our recent work 
D\v zamonja and Plebanek (2006) gives an easy
solution to the former of these problems, and gives some
partial information about the latter. The long term goal of this line of research
is to obtain a structure theory of Boolean algebras that support a finitely additive
strictly positive measure, along the lines of Maharam theorem which gives such
a structure theorem for measure algebras. 
\end{abstract}

\section{Introduction} This is an article in mathematics, specifically in set theory.
Set theory is a part of mathematical logic and it indeed has a dual role, that of
giving logical foundations to mathematics, and that of being a part of mathematics
itself. In many people's minds these two roles are rather distinct. In fact much of
the development of set theory in the twentieth century ran on two distinct
tracks. On the one hand, much effort was spent to develop systems
of set-theoretic axioms from which one could logically develop the known mathematics.
Hilbert's programme made it a priority to have such a system, and it had been
widely believed after Cantor (see e.g Cantor (1874)) developed set theory,
in the late
nineteenth century, that one should be able to have an axiomatic system for
mathematics by using
the notion of a set as given. The present understanding of this is rather
different, due mostly to the work of G\"odel (1931) in his famous
Incompleteness Theorems. Firstly he showed that for every consistent 
recursive system of axioms which includes the Peano 
Arithmetic there is a statement formalizable in the
theory itself which is independent (not provable nor unprovable)
in the theory, and secondly 
that such systems of axioms also cannot prove their own consistency.  
On the other hand, much
understanding was developed about the mathematical, rather than logical, consequences of
Cantor's work. This includes the idea of infinite sets of various sizes and the
understanding of what this means in terms of familiar objects, such as the sets of
reals. This article will mainly concentrate on that, mathematical, side of set theory,
but in fact the point is that the two sides of set theory are indivisible, 
as we wish to demonstrate by this article.

Most (but not all) of modern mathematics takes as a basis the axiom scheme known as
ZFC, Zermelo-Fraenkel axioms with Choice. It has been the case basically since the
1920s and there are many good reasons for this decision. Some then may view it as
a serious setback that it became known since Cohen (1963) that
not only that as we know from G\"odel (1931), there are some statements of
mathematics that are independent of ZFC, but that there are some ordinary
statements, notably the Continuum Hypothesis (`There is no infinite 
subset of the set ${\Bbf R}$ of real numbers which is not in a bijective correspondence
with either the set ${\Bbf N}$ of natural numbers or the set ${\Bbf R}$ itself')
that are independent. 
This was established by showing that
one cannot calculate just by using the axioms of ZFC the size $|{\Bbf R}|=|\PP({\Bbf N})|$
in terms of the Cantor's $\aleph${\em-hierarchy} of infinite cardinals:
$\aleph_0=|{\Bbf N}|$, the next infinite cardinal $\aleph_1$, $\aleph_2$, $\ldots$
the limit $\aleph_\omega$, the next $\aleph_{\omega+1}, \ldots$. Moreover,
it is consistent with these axioms (assuming they themselves are
consistent) that the value of $|\PP(\aleph_0)|$ is as large as desired, i.e. no upper
bound can be found just by arguing in ZFC.
Indeed, Cohen's result changed the subject
of set theory entirely, not only because of its logical significance but also because it
introduced a method for proving that various statements were independent of ZFC. This is
the method of forcing. Many interesting results have been obtained by applying this method.
Whilst this is exciting, it also created the feeling in the general mathematical
community that set-theorists are mostly concerned about things that cannot be done,
rather than the ones that can. It did not help that the forcing machinery
has developed to an incredible extent and even being able to read and verify some of
the proofs takes an enormous effort. There definitely was a period in which
set theory was considered far removed from the interests of the main stream mathematics.

I am very fortunate to belong to the generation of mathematicians who have now seen the
pendulum swing the other way. In recent years there have been a number of purely
mathematical results that have been obtained as a consequence of a fine set-theoretic
analysis of the problem. The final result often does not mention any
set-theoretic assumptions additional to ZFC, which is assumed throughout, yet the
proof relays deeply on an understanding of the set-theoretic limits of ZFC
and its possible universes. As an example, let us state a celebrated theorem of 
Shelah (1994):

\begin{theorem}(Shelah) If $|\PP(\aleph_0)|< \aleph_\omega$, then 
$|\PP(\aleph_\omega)|< \aleph_{\omega_4}$.
\end{theorem}

To appreciate the impact of the above theorem, contrast it with what we have already
said: it is not possible to bound $|\PP(\aleph_0)|$ by arguing in ZFC. However, if we are
in the situation to know that this value is less than $\aleph_\omega$ then we can 
put a definite bound on the size of $|\PP(\aleph_\omega)|$. 
This theorem as stated can be understood without any
prior knowledge of set theory. Yet, this statement is a culmination of at least twenty
years of concentrated effort by many set theorists, involving techniques such as large
cardinals, iterated forcing, elementary embeddings, and culminating by the seminal work
of Shelah (1994)
in which he invented the whole theory (`pcf', the theory of possible cofinalities)
to finally prove the theorem. 
This part of Shelah's work was what was cited when he
was awarded the prestigious Bolyai Prize. 

In this article we shall concentrate on the appearance of this phenomenon of the symbiosis
between the axiomatic and the mathematical in the context of measure theory. We shall 
describe the problem of characterizing Boolean algebras that carry a measure, and point
out the most well known instance of this problem. This is of course the von Neumann 
problem
(see von Neumann (1981)),
recently solved in the negative by Talagrand (2006). We shall then consider some other
instances of this problem and show some recent solutions that have been obtained in our
joint work D\v zamonja \& Plebanek (2006).

\section{Description of the problem}\label{known}
In order to make the article accessible we shall
commence with a quick review of the basic notions we shall use. A list of notational
conventions can be found at the end of the article.

A Boolean algebra $\algb$ 
is a structure consisting of a set with at least two distinct elements $0_{\algB}$ and
$1_{\algB}$, 
two binary operations, $\wedge$ and $\vee$, and a unary operation $-$, which obeys certain
rules known as the laws of Boolean algebras. A typical example of a Boolean algebra is
a family of subsets of a given set $A$, where $0_{\algB}=\emptyset$, $1_{\algB}=A$
and the operations $\wedge,\vee$ and $-$ are interpreted as $\cap,\cup$ and the complement
${}^c$ respectively. It follows from Stone
Representation Theorem (Stone (1936)) that every Boolean algebra is isomorphic as a structure
to some $\algB\subseteq \PP(A)$ for some $A$, so we shall only work with
such algebras. The basic laws of Boolean algebras are then interpreted as the
familiar commutativity, associativity and distributivity laws between $\cup$ and $\cap$,
and we also have that for any $a\in \algB$, $(a^{c})^c=a$. These operations induce
the familiar subset relation $\subseteq$, which acts as a relation of partial ordering
on $\algB$. Two elements $a,b$ of $\algB$ are said to be {\em incomparable} if
$a\cap b=\emptyset$.

Boolean algebras $\algB$ may also have properties additional to the ones
given by the basic laws. For example, we can
consider the {\em antichains}, which are subsets of $\algB$ consisting of pairwise
incomparable (i.e. disjoint)
elements. The condition that guarantees that all antichains in $\algB$ are
countable is called {\em the countable (anti)chain condition} and abbreviated as {\em ccc}.
An example of a Boolean algebra that satisfies this law is the family $\MM$ of 
the equivalence classes (mod. measure 0) of all Lebesgue
measurable subsets of the unit interval $[0,1]$. 
This is easily seen and is due to the
additivity properties of the Lebesgue measure $\lambda$. It is known 
(see Fremlin (1989))
that one can
choose the representatives $E^\bullet$ of the equivalence classes of measurable
sets $E$ so that 
$\emptyset^\bullet=\emptyset$ and the representative of $E^\bullet\cup (\cap)
F^\bullet$ is the union (intersection) of the corresponding representatives, for
all $E^\bullet, F^\bullet$.
(For this reason we omit ${}^\bullet$ in our notation). Another property of $\MM$ 
is that for every
sequence $\langle a_n:\,n<\omega\rangle$ in $\MM$ the union $\bigcup_{n<\omega}
a_n$ is an element of $\MM$, and it is the least upper bound of $\langle a_n:\,n<\omega\rangle$
with respect to $\subseteq$. Boolean algebras with this property are said to be 
$\sigma$-{\em complete}. In fact, the analogous completeness property remains
true for sequences indexed by any other ordinal but $\omega$, which can be proved by
using the ccc property along with the $\sigma$-completeness. This property is
called {\em completeness}, and since in this article we never deal with algebras that
are $\sigma$-complete without being complete, we shall simplify the notation and
refer to both concepts by the word `complete'. 

Notice that in the above
example $\lambda$ induces a function on $\MM$, which is again called $\lambda$, and that
this function satisfies the following, for all $a,b, a_n\in \MM$:
\begin{description}
\item{(i)} (strict positivity) $a\neq \emptyset\implies \lambda(a)>0$, $\lambda(\emptyset)=0$,
\item{(ii${}^-$)} (total finiteness) $\lambda(1_{\MM})<\infty$,
\item{(iii)} (additivity) if $a,b$ are disjoint then $\lambda(a\cup b)=\lambda(a)+\lambda(b)$,
moreover
\item{(iii${}^+$)} (countable additivity) if $\langle a_n:\,n<\omega\rangle$ are 
pairwise disjoint, then $\lambda(\bigcup_{n<\omega}a_n)=\Sigma_{n<\omega}\lambda(a_n)$.
\end{description}
A Boolean algebra $\algB$ which supports a functional $\lambda$ satisfying properties (i),
(ii)${}^{-}$
and (iii)${}^+$ above is called a {\em measure algebra} and
$\lambda$ is called a {\em strictly positive (s.p.) countably additive measure}. 
Measures which only satisfy properties (i), (ii)${}^{-}$ and (iii)
are called {\em s.p. finitely additive measures}. Since it is clear that
by multiplying by a constant we can obtain from $\lambda$ another countably
additive measure $\mu$ which satisfies $\mu(1_{\MM})=1$ (such measures are
called {\em probabilities}), we can replace the requirement (ii)${}^{-}$
in the definition of a measure algebra by the one
requiring the relevant measure to be a probability. In the sequel the word `measure'
will refer to finitely additive measures.

The general problem we shall discuss in this paper is the following:

\medskip

{\noindent {\bf Measure Recognition Problem}} MRP($\Phi)$ 
Given a Boolean algebra $\algB$ and a property
$\Phi$ of measures. How can we recognize by purely combinatorial means
if on $\algB$ one can define a strictly
positive measure with property $\Phi$? 

\section{Comments on the problem and known results}\label{knownres}

The most well known instance of the Measure Recognition Problem
was asked by von Neumann in 1937
(see von Neumann (1981)). He was interested in recognizing measure algebras (so the additional
property $\Phi$ in the problem description is the property of countable additivity
and we are dealing with MRP(countably additive)).
In addition to the completeness and the ccc property von Neumann isolated
another property which is always present in measure algebras, called weak distributivity,
and asked if these three properties together were sufficient for a Boolean algebra
to be a measure algebra. This famous problem was answered only very recently
by Talagrand (2006), and the answer is negative. In fact, the main
result of Talagrand (2006) answered,
also negatively,
the related well known Control Measure Problem, which asks if the existence of a
so called continuos submeasure on Boolean algebra implies the existence of a measure.
The step from this solution to the solution of von Neumann's problem then followed
by known work, as explained in Talagrand's paper. A 
(strictly positive) continuos submeasure is a finite
nonnegative function $\mu$ on a Boolean algebra $\algB$, vanishing only at $0_{\algB}$,
having the `submeasure' property that $\mu(a\cup b)\le \mu(a)+\mu(b)$ for all
$a,b\in \algB$, and the `continuity' property that for any sequence 
$\langle a_n:\,n<\omega\rangle$ of elements of $\algB$ satisfying $a_{n+1}\supseteq a_n$
for all $n$
and $\bigcap_{n<\omega} a_n=\emptyset$,
we have $\lim_{n}\mu(a_n)=0$. The result of Talagrand (2006) is probably
the most exciting recent result to come out of a whole variety of mathematical subjects,
particularly measure theory and set theory. 

The notion of a continuous submeasure comes
from the work of Maharam (1947) who observed that a necessary condition for a Boolean
algebra to be a measure algebra is to be metrizable, and showed how one can define
a continuos submeasure from the appropriate metric. Following this work an important
question became to recognize Boolean algebras that support continuous submeasures,
the so called Maharam's algebras.
Recently it was proved in Balcar \textit{et al.} (2005) and Veli\v ckovi\'c (2005)
that under a certain set-theoretic axiom known as the P-ideal
dichotomy, any ccc weakly distributive complete Boolean algebra is a Maharam algebra.
Using their work as a starting point Todor\v cevi\'c (2004) obtained the following
interesting characterization:

\begin{theorem} (Todor\v cevi\'c) A complete Boolean algebra carries a strictly positive
continuous submeasure if and only if it is weakly
distributive and satisfies the $\sigma$-finite chain condition.
\end{theorem}

The chain condition here means that the algebra can be written as a countable union of
subsets each of which only has finite antichains. Todor\v cevi\'c's result is clearly
a theorem of ZFC, but it was obtained as a consequence of a set-theoretic analysis
of the situation under the P-ideal dichotomy, and the methods introduced in
Maharam (1947). Likewise, Talagrand's result is
a theorem of ZFC but even the formulation by Maharam
of the Control Measure Problem was a consequence of her analysis of the behaviour
of the von Neumann's problem under the set-theoretic assumption of the existence of a
Souslin tree. Both of these results illustrate the point that this article makes, about
the close connection between a set-theoretic analysis of a problem and its solution in
ZFC, as well as the long-term view that one has had to take in understanding various specific
instances of this connection.

Next let us take a look at another special case of MRP, namely MRP($\emptyset$),
the situation when we do not require any special property $\Phi$. In this case there is
a combinatorial criterion due to Kelley (1959). It involves the notion of the 
{\em intersection number}
of a family $\FF$ of sets, which is defined to be the supremum of all $\alpha$ such that 
for every finite subsequence $\bar{a}$ of elements of $\FF$ (with possible repetitions of
elements), there is a subsequence $\bar{b}$ of length $\lg(\bar{b})$
at least $\alpha\cdot\lg(\bar{a})$, such that the intersection of all elements of $\bar{b}$
is non-empty. Kelley's criterion is then:

\begin{theorem} (Kelley) A Boolean algebra $\algB$ carries a strictly positive
(finitely additive) measure if and only if $\algB\setminus\{\emptyset\}$ can be written
as a countable union of families each of which has positive intersection number.
\end{theorem}

One may wonder how strong the condition in the Kelley's criterion is. It clearly 
implies the ccc, but Gaifman (1964) proved that there are ccc algebras that do not
satisfy Kelley's criterion.
It is also known that under the set-theoretic statement known as Martin's Axiom MA
and the negation of the continuum hypothesis CH, all ccc Boolean algebras of size
$<2^{\aleph_0}$ do satisfy Kelley's condition and in fact some stronger
conditions (see Fremlin (1984)).
This together with the example of the P-ideal dichotomy results quoted above 
demonstrates that
there are some mathematical axioms which make it easier for a Boolean algebra to have
certain measure-theoretic properties. We shall see another example of this behaviour in
\S\ref{new}. Once the notion of the intersection number is known, the proof of Kelley's
theorem follows rather readily by the well known facts from functional analysis.
It seems to be a tendency in this subject that
positive theorems once formulated properly, have proofs that are much less involved
than the proofs of the negative ones. We shall illustrate this in Theorem \ref{baanonatomic},
where we shall give a proof of one such positive theorem, while we may note that
most of the negative statements mentioned here (for example
the construction from Talagrand (2006)) have proofs that are very much out of the scope
of this paper.

Special cases of property $\Phi$ which are particularly interesting
are the notions of nonatomicity and separability. To motivate the definitions
we shall go back to the case of countably additive measures. An example of such a
measure is the familiar Lebesgue measure $\lambda$ on the unit interval. This measure naturally
leads to a measure on the Tychonoff product $[0,1]^\kappa$ for any cardinal $\kappa$,
denoted by $\lambda_\kappa$. A remarkable theorem of Maharam (1942) shows that 
the structure of measure algebras (so supporting a countably additive measure) is totally 
determined by these examples. Namely, any measure algebra can be decomposed into a
countable union of subsets, each of which is either an atom, or a Boolean algebra
isomorphic to the measure algebra $\algB_\kappa$
of some $[0,1]^\kappa$ under $\lambda_\kappa$, where all $\kappa$ are infinite cardinals.
Here we use the notion of an {\em atom} of a Boolean algebra
$\algB$, which is an element $a\neq \emptyset$
such that only $b\subseteq a$ in $\algB$
are $\emptyset$ and $a$. Such a structure theorem does not
exist for Boolean algebras that simply support a strictly positive finitely additive
measure. Maharam (1981) writes: `It would be very interesting to have a structure 
theory for finitely additive measures, the structures of which can be much more
complicated (than those of measure algebras).'

Taking a long term goal to obtain such a structure theory for finitely additive measures,
we may ask ourselves what the first step would be. In the case of countably additive
measures that was to determine the relevant building blocks, namely the algebras $\BB_\kappa$.
Each of these blocks has two important properties. The first one is that the
measure is {\em nonatomic} which means that for every
$\eps>0$ there is a finite partition of the algebra into elements of measure $<\eps$. 
For this reason we wish to have a combinatorial criterion for recognizing when a
Boolean algebra has a finitely additive measure which is nonatomic, which explains
why we believe the case of MRP(nonatomic) is an important special
case. In our recent work D\v zamonja \& Plebanek (2006) we obtained a simple
solution to this problem, which
will be presented in \S\ref{new}.

Another property of measure algebras $\BB_\kappa$ for
$\kappa\ge\aleph_0$ is that they can be understood as metric spaces of
density $\kappa$. Namely one introduces a metric $d_\kappa$ derived from $\lambda_\kappa$
by declaring $d_\kappa(a,b)=\lambda_\kappa(a\Delta b)$, and this metric has
the required property. By Maharam's theorem mentioned above, every measure algebra $\algB$
has a (unique) decomposition involving algebras $\algB_\kappa$ for some $\kappa$, and the
supremum of all $\kappa$ involved in this composition is called {\em Maharam's dimension}
or {\em type}
of $\algB$. Hence it would be of interest to have a similar notion for Boolean algebras
which simply support a finitely additive strictly positive measure. The notion of type
can be defined similarly to the above because a finitely additive strictly positive
measure on $\algB$ will already induce a metric on $\algB$, and we can define
the {\em type} of $\algB$ to be the supremum of all densities of metric spaces obtained
from $\algB$ by using all possible strictly positive measures on $\algB$. In the same vein,
for a fixed measure $\mu$ on $\algB$ we may consider the density of the induced metric
space. In particular we say that $\mu$ is {\em separable} if the induced metric space
is separable in the topological sense (i.e. it has a countable dense subset).
Recognizing Boolean algebras that have such a strictly positive measure translates in the context
of MRP into MRP(separable). Clearly, a similar notion can be defined for
any fixed possible densities of a metric space, but for the moment we still do not
know how to solve MRP(separable). We shall show some partial solutions in
\S\ref{new}.

\section{Some recent results}\label{new}
In this section we shall take for granted all notions defined in previous sections.
Here we concentrate on MRP with two specific values of $\Phi$ in mind,
MRP(nonatomic) and MRP(separable). All otherwise unattributed theorems are from
D\v zamonja \& Plebanek (2006). This work in particular solves
the problem MRP(nonatomic) mentioned above, using a rather simple argument. We present
the argument here:

\begin{theorem}\label{baanonatomic}  A Boolean algebra $\algB$ carries
a strictly positive nonatomic measure 
if and only if there is a decomposition $\algB\setminus\{0\}=
\bigcup_{n<\omega} \algB_n$, where for each $n$ we have

\begin{itemize}
\item[(i)] $\algB_n\sub \algB_{n+1}$;
\item[(ii)] $\intn(\algB_n)\ge 2^{-n}$;
\item[(iii)] if $a\in \algB_n$
then there are disjoint $b,c\in \algB_{n+1}$
with $b\cup c\le a$. 
\end{itemize}
\end{theorem}

The forward direction of this theorem is easy modulo known facts. Namely, one uses
the well known Stone duality theory
between Boolean algebras and compact zerodimensional topological spaces, Stone (1936) and
its application to measures to transfer the
problem into the setting of Radon measures on compact spaces. The conclusion then
follows by another well known theorem, Maharam (1942). Details are not practical to
explain here.
We shall however
sketch the proof of the backward direction of the theorem, assuming the following lemma
which appears as part of Kelley's proof in
Kelley (1964) and os taken in this form
from Fremlin (2002), Proposition 391 I. The notation 
$\intn(A)$ stands for the intersection number of the family $A$. 

\begin{lemma}\label{int} (Kelley)
Let $\alg$ be a Boolean algebra and $A\subseteq \alg\setminus\{0\}$
nonempty. Then 
\[
\intn(A)=\max_{\nu}\inf_{a\in A}\nu(A),
\]
where the maximum is taken over all  probability (finitely additive) measures on $\alg$.
\end{lemma} 

\begin{proof} 
If there is a decomposition of $\algB$ satisfying (i)--(iii) then
by Fact \ref{int} for each $n$
we can define a probability measure $\mu_n$ on $\algB$ such that for all $b\in\algB_n$
we have $\mu_n(b)\ge 2^{-n}$. 
We let $\mu$ be any cluster point of the sequence
$\langle\mu_n:\,n<\omega\rangle$. It is easily seen that $\mu$ is a probability measure on $\algB$.
Let us show that $\mu$ is strictly positive.
By induction on $n$ it easily follows that for all $a\in \BB_n$ and $k\ge n$ 
there are $2^{k-n}$ pairwise disjoint sets in $\BB_{k}$ contained in $a$.
Hence for such $n,k$ we have $\mu_k(a)\ge 2^{k-n}\cdot 2^{-k}=2^{-n}$. Consequently,
$\mu(a)\ge 2^{-n}>0$. 

Suppose now that $\varepsilon>0$ and $a\in \algB\setminus\{0\}$ are given. 
Let $n\le k$ be large enough so that $1-2^{-n} +2^{-k}<\varepsilon$
and $a\in \BB_n$. Let $b_0,\ldots b_{2^{k-n}-1}$
be disjoint elements of $\BB_{k}$ contained in $a$, which exist as shown in the
previous paragraph. Then $\mu(b_i)\ge 2^{-k}$ for all $i< {2^{k-n}-1}$,
again by the argument in the previous paragraph, and hence $\mu(b_0)\le
1-[(2^{k-n}-1)\cdot 2^{-k}]= 1-2^{n}+2^{-k}<\varepsilon$.
This shows that every nonzero element of $\algB$ has a subset of positive measure
less than $<\varepsilon$, which easily implies nonatomicity.
\end{proof}

An essential difference between this criterion and other combinatorial criteria we
mentioned above (Maharam, Kelley, Todor\v cevi\'c), is that the decomposition of the Boolean
algebra involves an interaction between the countably many pieces involved. Next, one may 
wonder how strong this criterion is. Obviously, a necessary condition on the Boolean
algebra to support a nonatomic strictly positive measure is first of all that it
supports any strictly positive measure, so Kelley's criterion, and secondly that
the algebra itself is {\em atomless} which means that it does not have any atoms in
the sense defined in \S\ref{knownres}. D\v zamonja \& Plebanek (2006) give an
example showing that in general these conditions are not sufficient to
show that the algebra supports a nonatomic strictly positive measure. 
As it is known that under the axioms MA without CH
various properties of Boolean algebras of size smaller than $2^{\aleph_0}$ tend
to be equivalent, the following theorem is perhaps not surprising:

\begin{theorem}\label{chaincond} Assume MA and the negation of CH. Then for atomless Boolean 
algebras $\BB$ of size $<2^{\aleph_0}$, the following are equivalent
\begin{description}
\item{(i)} $\BB$ is ccc, and
\item{(ii)} $\BB$ satisfies the condition from Theorem \ref{baanonatomic}.
\end{description}
\end{theorem}

What may be surprising is the way that the theorem is proved. We give
an informal sketch, for which we need several notions. This will also 
be a good opportunity to introduce an auxiliary property of Boolean algebras
which was central in D\v zamonja \& Plebanek (2006), the so called
approximability.

A {\em forcing notion} is
a partially ordered set $\mathbf P$ with the least element. Two elements of 
$\mathbf P$ are {\em incompatible} if there is no element of $\mathbf P$ that is
larger than both of them. We say that $\mathbf P$ is {\em ccc} if there is no 
uncountable family of pairwise incompatible elements (note that this is a different
notion than that of ccc in Boolean algebras). A subset of $\mathbf P$ is a {\em filter}
if it is directed and downward closed as a partial order. A subset $D$ of $\mathbf P$ is
{\em dense} if for every $p\in {\mathbf P}$ there is $q\in D$ with $p\le q$.
 MA states that for any ccc forcing $\mathbf P$
and any family $\FF$ of $<2^{\aleph_0}$ dense subsets of $\mathbf P$, there is a filter
of $\mathbf P$ which intersects each of the dense sets in the family. In applications
one may start with a goal of constructing some object, usually of size $\aleph_1<2^{\aleph_0}$
(which is why the negation of CH is assumed as well), and $\mathbf P$ consisting of
(usually) finite pieces of the desired object, ordered so that a stronger element of
$\mathbf P$ gives more information about the object than a weaker one. Then one need to
formulate a family $\FF$ of dense sets which represent various requirements
on the object constructed, and the $\FF$-generic filter guaranteed to exist by MA will
in some natural way give rise to the object desired. 

In the context of Theorem \ref{chaincond}, in the nontrivial direction from
(i) to (ii) and in view of characterization from theorem \ref{baanonatomic}, it would
be natural to start with a Boolean algebra $\algB$ which satisfies ccc and has size $<2^{\aleph_0}$,
and to formulate a ccc forcing notion which would force a nonatomic 
strictly positive measure on $\algB$. One would need to formulate a dense set corresponding
to the requirement that $a$ has positive measure, for every nonzero $a$ in $\algB$.
This is maybe possible, but this is not how our proof goes. Also, the known proof
that every ccc Boolean algebra of size $<2^{\aleph_0}$ is
separable and hence satisfies Kelley's condition
and carries a measure (see Fremlin (1984)), is not useful since it gives an atomic measure.
We worked instead with the
notion of approximability. A Boolean algebra $\algB$ is {\em approximable} if
there is a sequence $\langle \mu_n:\,n<\omega\rangle$ of probability
measures on $\algB$ such that for every $a\neq 0$ in $\algB$ there is $n$ such that
$\mu_n(O)>1/2$. This notion was studied by Talagrand (1980) and 
M\"agerl \& Namioka (1980) with the idea of characterizing Boolean algebras that
support a strictly positive separable measure. We shall discuss this in a
moment, but for now let us finish sketching the proof of Theorem \ref{chaincond}.
The proof is to start with a Boolean algebra that satisfies Kelley's condition
(in fact a weaker condition than that) and formulate a ccc forcing notion whose
generic filter gives a sequence $\langle \mu_n:\,n<\omega\rangle$ which demonstrates that
the algebra is approximable. The size of the family $\FF$ of
sets such that we need to produce an $\FF$-generic filter is equal to the
size of $\algB$, so if MA and not CH holds we can deal with any ccc algebra 
of size $<2^{\aleph_0}$. Having the sequence $\langle \mu_n:\,n<\omega\rangle$
we can define the weighted sum $\Sigma_{n<\omega}\mu_n/2^{n+1}$, and the proof is done in
such a way that if the algebra is atomless then this measure in nonatomic. 

Having sketched the proof of Theorem \ref{chaincond} we can discuss the next item on
our list, namely the notion of separable measures defined above. Talagrand (1980)
attacks MRP(separable) by defining a list of properties of decreasing strength one
of which is `supporting a strictly positive separable measure' (SM). He shows that
the property of $\algB$ being {\em $\sigma$-centred} (so the algebra can be written as a
countable union of subfamilies each of which has the property that for all finite
subsets $J$ we have $\bigcap J\neq \emptyset$, provided
$\emptyset\notin J$) is strictly
stronger than SM, while SM implies approximability. He also shows that under the
assumption of CH, property SM is strictly stronger than approximability. 
There is a possibility left open by this theorem,
which is that approximability might under some suitable axioms actually characterize
SM. D\v zamonja \& Plebanek (2006) show that this is not the case:

\begin{theorem}\label{Grz} There is a Boolean algebra which is approximable but does not
support a strictly positive separable measure.
\end{theorem}

Proof of this theorem involves a combinatorial construction the details of which
are out of the scope of this paper. A major building block is a zerodimensional
topological space without isolated points
with the curios property that it has a countable dense set
$D$ such that for every sequence $\langle F_n:\,n<\omega\rangle$ of closed sets
whose union is disjoint from $D$, that union is nowhere dense. Such a space was
first constructed in ZFC by Simon (2002), although it was known before that such spaces can
exist under various additional axioms of set theory.

Just as Talagrand (2006) by answering negatively the problem of von Neumann
brings us back to square one regarding MRP(countably additive), so does
Theorem \ref{Grz} brings us back to the beginning regarding MRP(separable).
It would have been very nice if at least consistently there would be an equivalence
between approximability and separability, since there is a combinatorial characterization
of approximability due to M\"agerl \& Namioka (1980):

\begin{theorem}\label{baapproximable} (M\"agerl-Namioka) 
A Boolean algebra $\algB$ is approximable
if and only if for every $\varepsilon>0$ (equivalently: for some $\eps\in (0,1))$
there is a decomposition $\algB\setminus\{0\}=
\bigcup_{n<\omega} \algB^\varepsilon_n$, where for each $n$ we have
$\intn(\algB_n^\varepsilon)\ge 1-\varepsilon$.
\end{theorem}

In the absence of such an equivalence we have to rethink what a possible characterization
may look like. At this point we may say that it is likely that any characterization
will involve not only a decomposition into countably many pieces each of which has
some given property, in the style of all but one characterization given here. Instead,
the desired theorem will probably have to take into account the interaction between these
various pieces, like in Theorem \ref{baanonatomic}. The reason for this is
a result by Dow \& Steprans (1993). Namely,
one of the strongest possible characterizations involving only the information about each 
of the countably many pieces (given that we know that $\sigma$-centredness is too strong),
would use the notion of {\em $\sigma-n$-linkedness}. For a given $n$, a family $A$ of
nonempty sets is $n$-linked if for any subset $J$ of $A$ of size $\le n$, we have
$\bigcap J\neq\emptyset$. A Boolean algebra is {\em $\sigma-n$-linked} if it can be
written as a countable union of sets each of which is either $\{\emptyset\}$ or
$n$-linked.
Dow \& Steprans (1993) proved that the measure algebra of $[0,1]^\kappa$ is
$\sigma-n$-linked for every $n$ if and only if $\kappa\le 2^{\aleph_0}$. This shows that
the notion of $\sigma-n$-linkedness does not make a sufficient distinction between
separable measure algebras and those of type $\le 2^{\aleph_0}$, and it is easily
seen using Stone's duality (from Stone (1936))
that the measure algebra of $[0,1]^\kappa$ does not support a strictly positive
separable measure for $\kappa>\aleph_0$.

\section{Conclusion} The article shows on the specific example of the Measure Recognition
Problem the phenomenon of the utility of a multidisciplinary mathematical approach to
a single mathematical problem, in particular the value of a set-theoretic analysis
of the problem at hand. 
We have shown how MRP, which is a problem about Boolean algebras
that was asked first at least in 1937 if not earlier, has had an impact of measure theory,
combinatorics and set theory, and in turn that each these subjects contributed
to a better understanding of MRP. Specifically, we discussed the set-theoretic
insights and emphasized the point of using
set-theoretic tools and being able to obtain 
a purely measure-theoretic or algebraic final result.
Through historical remarks presented in the first
parts of the paper the reader can see how various results have followed the current state of
knowledge, both in set theory and measure theory. In the final section we showed 
some new results and underlined the future directions of research in the subject.
One of the goals is a structure theory for Boolean algebras that support a finitely
additive strictly positive measure.

\section{Notation} The {\em power set} of a given set $A$ is the set of all subsets of
$A$ and is denoted by $\PP(A)$.

The set-theoretic notation $\{ E_n:\,n<\omega\}$ is used for
what some authors denote by $\{ E_n:\,n\in {\Bbf N}\}$ or (after a reenumeration)
$\{ E_n:\,1\le n <\infty\}$. Our notation for sequences is
$\langle a_n:\,n<\omega\rangle$ in place of $(a_n)_n$. A similar convention is used for
$\Sigma_{n<\omega}$.

Many notions were introduced in the text, where at first appearance they were italicized.
\medskip
\medskip

\begin{acknowledgements} The author thanks EPSRC for their support through an Advanced Fellowship
in Mathematics and the British Council for support through an Alliance Grant for years 2005
and 2006.
\end{acknowledgements}

\begin{thedemobiblio}{}

\item Balcar, B., Jech T. \& Paz\'ak T. 2005. Complete ccc boolean algebras,
the order sequential topology, and a problem of von Neumann. \textit{Bull.
London Math. Soc.} \textbf{137}, 6, 885--898.

\item Cantor, G. 1874. \"Uber eine Eigenschaft des Inbergriffs aller reellen algebraischen
Zahlen. \textit{J. f. Math.} \textbf{77}, 258--262.

\item Cohen, P. 1963. The independence of the continuum hypothesis. \textit{Proc. Natl. Acad. Sci
USA} \textbf{50}, 1143--1148.

\item Dow, A. \& Steprans, J. 1993. The $\sigma$-linkedness of measure algebra.
\textit{Can. Math. Bull.} \textbf{37}, 1, 42--45.

\item D\v zamonja, M. \& Plebanek, G. 2006. Measures on Boolean algebras. Preprint. University
of East Anglia and Wroc\l aw University.

\item Fremlin, D.H. 1984. \textit{Consequences of Martin's Axiom}. Cambridge:
Cambridge University Press.

\item Fremlin, D.H. 1989. Measure algebras. In \textit{Handbook of
Boolean algebras} \textbf{III} (ed. J. D. Monk) 877--980. Amsterdam: North--Holland. 

\item Fremlin, D.H. 2002. \textit{Measure Theory (Measure
Algebras)}, {\textbf 3}, Colchester: Torres Fremlin.

\item Gaifman, H. 1964. Concerning measures on Boolean algebras. \textit{Pacific J. Math.},
\textbf{14}, 61--73.

\item G\"odel, K.  1931. \"Uber formal unentscheidbare S\"atze der Principia Mathematica und
verwandter Systeme. I. \textit{Monatsh. Math. Phys.} \textbf{38}, 173--198.

\item Kelley, J. L. 1959. Measures on Boolean algebras. \textit{Pacific J. Math.} \textbf{9},
1161--1177.

\item M\"agerl, G. \& Namioka I. 1980.
Intersection numbers and weak${}^\ast$ separability of spaces of measures.
\textit{Math. Ann.} {\textbf 249}, 3, 273--279.

\item Maharam, D. 1942. On homogeneous measure algebras. \textit{Proc. Nat. Acad. Sci. USA}
\textbf{28}, 108--111.

\item Maharam, D. 1947. An algebraic characterization of measure algebras. \textit{Ann. of Math.}
\textbf{48}, 2, 154--167.

\item Maharam, D. 1981. Commentary on Problem 163. In \textit{The Scottish Book}
(ed. R. D. Mauldin) 241--243. Boston: Birkha\"user.

\item Shelah, S. 1994. \textit{Cardinal Arithmetic}. Oxford: Oxford University Press.

\item Simon, P. 2002. A countable dense-in-itself dense P-set.
\textit{Topol. Appl.} \textbf{123}, 193--198.

\item Stone, M.H. 1936. The theory of representations for Boolean algebras. \textit{Trans.
Amer. Math. Soc.} \textbf{40}, 37--111.

\item Talagrand, M., 1980. S\'eparabilit\'e vague dans l'espace des mesures sur un
compact. \textit{Israel J. Math.} \textbf{37}, 1-2, 171--80.

\item Talagrand, M. 2006. Maharam's Problem. Preprint. University of Paris and Ohio
University.

\item Todor\v cevi\'c, S. 2004. A problem of von Neumann and Maharam about
algebras supporting continuos submeasures. \textit{Fund. Math.}
\textbf{183}, 2, 169--183.

\item Veli\v ckovi\'c, B. 2005 CCC forcings and splitting reals. \textit{
Israel J. Math.} \textbf{147}, 209--220.

\item von Neumann, J. 1981. Problem 163. In \textit{The Scottish Book} (ed. R. D. Mauldin) pg.~240.
Boston: Birkha\"user.

\end{thedemobiblio}
\end{document}